\documentclass[12pt,reqno]{amsart}
\usepackage{mathrsfs}
\usepackage{amsfonts}
\usepackage{amsmath}
\numberwithin{equation}{section}

\newtheorem{theorem}{Theorem}[section]
\newtheorem{lemma}{Lemma}[section]
\newtheorem{proposition}{Proposition}[section]
\newtheorem{definition}{Definition}[section]
\newtheorem{remark}{Remark}[section]

\title[\ Global   solution  to  nonlinear  Dirac  equation]
      {Global   solution  to  nonlinear  Dirac  equation for Gross-Neveu model in $1+1$ dimensions}

\author[ ]{Yongqian Zhang and Qin Zhao}
\keywords{Nonlinear Dirac equation, Gross-Neveu model, global strong solution, Bony type functional, Glimm type functional  \\  ( AMS subject classification. Primary: 35Q41 ; Secondary: 35L60, 35Q40) }
\email{yongqianz@fudan.edu.cn }
\email{qinzhao11@fudan.edu.cn}

\begin{document}
\maketitle

\centerline{\scshape }

{\footnotesize \centerline{School of Mathematical Sciences }
   \centerline{ Fudan University, Shanghai 200433, P.R.China}
} 


\begin{abstract}
This paper studies a class of nonlinear  Dirac equations with cubic terms in $R^{1+1}$, which  include the equations for the massive Thirring model and the massive Gross-Neveu model. Under the assumption that the initial data has bounded $L^2$ norm, the global existence and the uniqueness of the strong solution in $C([0,\infty),L^2(R^1))$ are proved.
\end{abstract}

\section{Introduction}

We consider the nonlinear  Dirac equations
\begin{equation}\label{eq-dirac}
\left\{ \begin{array}{l} i(u_t+u_x)=-mv +N_1(u,v), \\i( v_t-v_x)=-mu+N_2(u,v),
\end{array}
\right.
\end{equation}
with initial data
\begin{equation}\label{eq-dirac-initialv}
(u, v)|_{t=0}=(u_0(x), v_0(x))
\end{equation}
where $(t,x)\in R^2$, $(u,v)\in \mathbf{C}^2$, $m\ge0$. The nonlinear terms  takes the following form
\begin{equation}\label{eq-nonlinearstruc} N_1=\partial_{\overline{u}}W(u,v), \quad N_2=\partial_{\overline{v}}W(u,v)
\end{equation}
with
\[ W(u,v)=\alpha |u|^2|v|^2+\beta (\overline{u}v+u\overline{v})^2,\]
where $\alpha,\beta\in R^1$ and $\overline{u}, \overline{v}$ are complex conjugate of $u$ and $v$.
  (\ref{eq-dirac}) is called  Thirring equation for $\alpha=1$ and $\beta=0$, while it is called Gross-Neveu equation for $\alpha=0$ and $\beta=1/4$; see for instant \cite{thirring} and \cite{gross-neveu}, \cite{pelinovsky}.

 There are a number of works on the local and global well-posedness of Cauchy problem  for nonlinear  Dirac equations with cubic or higher-order nonlinear terms in different Sobolev spaces, see for examples \cite{bachelot}, \cite{bachelot-2}, \cite{bournaveas-zouraris}, \cite{cacciafesta1}, \cite{cacciafesta2}, \cite{candy},  \cite{deldado},  \cite{dias-figueira1}, \cite{dias-figueira2}, \cite{escobedo}, \cite{Esteban-Lewin-Sere}, \cite{huh}, \cite{huh2}, \cite{machihara}, \cite{pelinovsky}, \cite{selberg}, \cite{zhang} and the references therein. The survey of the well-posedness in the nonlinear Dirac equation in one dimension is given in \cite{pelinovsky}. Recently, the global existence of solutions in $L^2$ for Thirring model (i.e. $\beta=0$) has been established by Candy  in \cite{candy}, while the  global existence in $C([0,+\infty); H^s(R^1))\cap C^1([0,+\infty); H^{s-1}(R^1))$ for $s>1/2$ has been obtained by Huh in \cite{huh2} for (\ref{eq-dirac}).  For the bounded initial data with small $L^2$ norm, the existence and the uniqueness of global strong solution in $C([0,\infty); L^2(R^1))\cap L^{\infty}([0,\infty)\times R^1)$ for massless Gross-Neveu model (i.e. $m=0$) has been given in \cite{zhang}. To our knowledge, the global existence of  solution in $L^2$  for the Dirac equation of Massive Gross-Neveu model is still open \cite{pelinovsky}.
This paper is to devoted to study this problem. We consider a class of cubic nonlinear Dirac equations in one dimension, which include the equation for Massive Gross-Neveu model as an example. For  initial data with bounded $L^2$ norm, we will find the global solution $(u,v)$ which solves (\ref{eq-dirac}) and (\ref{eq-dirac-initialv}) in following sense.
  \begin{definition}\label{def-weaksolution}
  A pair of functions $(u,v)\in C([0,\infty); L^2(R^1))$ is called a strong solution to (\ref{eq-dirac}) and (\ref{eq-dirac-initialv}) on $R^1\times [0,T]$ if there exits a sequence of classical solutions $(u^{(n)}, v^{(n)}) $ to (\ref{eq-dirac}) $R^1\times [0,T]$ such that
  \[ (u^{(n)},v^{(n)}) \to (u,v), \quad \mbox{strongly in $L^2(R^1\times [0,T])$,}\]
  \[ (u^{(n)},v^{(n)})\big|_{t=0} \to (u_0,v_0), \quad \mbox{strongly in $L^2(R^1)$}\]
  as $n\to \infty$.
  \end{definition}
  \begin{remark}
  For the sequence of smooth solutions $(u^{(n)},v^{(n)})$ given in Definition \ref{def-weaksolution},   we can show as in the proof of Proposition \ref{prop-convergence2} in Section 3 that
  \[ u^{(n)}v^{(n)} \to uv \quad \mbox{in} \, \,  L^2_{loc}(\overline{R^1\times R^1_+}),\]
  which leads to
  \[ N_j(u^{(n)},v^{(n)}) \to N_j(u,v), \, \, j=1,2, \quad \mbox{in} \, \, L^1_{loc}(\overline{R^1\times R^1_+}).\]
  Therefore $(u,v)$ is a weak solution to (\ref{eq-dirac}) and (\ref{eq-dirac-initialv}), that is, a strong solution $(u,v)$ is also a distributional  solution in the standard sense.
  \end{remark}
  To construct the global strong solution for  (\ref{eq-dirac}) and (\ref{eq-dirac-initialv}), we approximate the initial data (\ref{eq-dirac-initialv}) by  a sequence of smooth functions. Then,
  by Huh's recent work \cite{huh2} on the global wellposedness  in $C([0,T],H^1(R^1))$ for (\ref{eq-dirac}), we can have a sequence of global smooth solutions for smooth data for (\ref{eq-dirac}). Our aims in this paper is to prove the convergence of the smooth solutions and to prove the uniqueness of the strong solution. To this end, we first establish some estimates for smooth solutions $(u,v)$. We recall a  result on the pointwise estimates for $|u|^2$ and $|v|^2$, see Lemma \ref{lemma-PointwiseEstimate}; then write down the equations (\ref{eq-dirac1}) for $|u|^2$ and $|v|^2$  and the equations (\ref{eq-dirac-difference1}) for the difference of smooth solutions, which play the important roles in our analysis. We remark that (\ref{eq-dirac1}) has been used by Huh \cite{huh2} and Zhang \cite{zhang} to obtain the $L^{\infty}$ bounds on the solutions for massive or massless Gross-Neveu equations and that $L^{\infty}$ bounds on the solutions for Thirring model has been obtained by Delgado \cite{deldado} in similar way. In contrast to \cite{deldado}, \cite{huh2} and \cite{zhang}, $L^{\infty}$ bounds on the $L^2$-strong solutions seems unavailable for (\ref{eq-dirac}) since there is no $L^{\infty}$ control on the initial data here. Therefore  we present two new functionals  to get the control on the nonlinear terms  in (\ref{eq-dirac1})   and  (\ref{eq-dirac-difference1}) for smooth solutions.
 The Bony type functional $\displaystyle Q_0(t)=\mathop{\int\int}_{a-t_0+t<x<y<b+t_0-t} |u(x,t)|^2|v(y,t)|^2dxdy$ is introduced to get the uniform $L^2_{loc}$ bounds on $uv$ for smooth solution $(u,v)$, while  the Glimm type functional $L_1(t)+KQ_1(t)$ is introduced to derive the local $L^2$ estimates on the difference between these smooth solutions, see Definitions \ref{def-bonyfunctional} and \ref{def-bony-difference}. These estimates enable us to get the local convergence in $L^2_{loc}$ of the sequence of smooth solutions  and the convergence in $L^1_{loc}$ of nonlinear terms corresponding to the smooth solutions, which give the local strong solution. Then we make the use of Lemma \ref{lemma-PointwiseEstimate} to extend the local solution globally. The uniqueness of solution could be proved in the same way.
 We remark that the functionals $Q(t)$ and $Q_1(t)$ are analogous to the Glimm's interactionl potential \cite{glimm} and the Bony functional \cite{bony}. The  Glimm interaction potential was first used by Glimm \cite{glimm} to establish the global existence of the  solutions with bounded variations for the system of conservation laws, then was used to study the general conservation laws; see for instance \cite{glimm},  \cite{bressan}, \cite{dafermos}, \cite{tartar} and the references therein. Bony type functional was used to get global solution in $L^1$ for Boltzmann equations and to study the wave maps in $R^{1+1}$, see for instance  \cite{bony},   \cite{ha-tzavaras}and \cite{zhou} and the references therein.

The main result is stated as follows.
\begin{theorem}\label{thm-existence}
For $(u_0,v_0)\in L^2(R^1)$, the problem (\ref{eq-dirac}-\ref{eq-dirac-initialv}) admits a unique global strong solution $(u,v)\in C([0,\infty);  L^2(R^1))$; moreover $|u||v|\in L^2([0,T]\times R^1)$ for any $T>0$.
\end{theorem}

The remaining of the paper is organized as follows. In Section 2 we recall Huh's result on the global well-posedness in $C([0,\infty); H^s(R^1))$ for $s>1/2$ for (\ref{eq-dirac}) at first.  Due to this result, we consider the smooth solutions at first. Then the Bony type functionals $Q_0(t)$ and $Q_1(t)$ are introduced for smooth solutions to deal with nonlinear terms in (\ref{eq-dirac1}) and (\ref{eq-dirac-difference1}-\ref{eq-dirac-difference21}), see Definitions \ref{def-bonyfunctional} and \ref{def-bony-difference}. With these functionals, we obtain the bounds on $L^2$ for smooth solution and the estimates on difference between smooth solutions on the intervals with small lengthes. In Section 3 the proof of the main result is given by  the estimates in Section 2.

\section{Estimates on the global classical solutions }\label{section-estimates-solu}

\subsection{Estimates on the classical solutions along the characteristics}

In the recent paper \cite{huh2}, Huh  established the global existence of the solution in $C([0,+\infty); H^s(R^1))\cap C^1([0,+\infty);H^{s-1}(R^1))$ to (\ref{eq-dirac}) for $s>\frac{1}{2}$. The main result in \cite{huh2} can be stated as follows.
\begin{theorem}\label{thm-huh}
For $(u_0,v_0)\in H^s(R^1)$ with $s>1/2$, (\ref{eq-dirac}-\ref{eq-dirac-initialv}) has a unique global solution $(u,v)\in C([0,\infty); H^s(R^1))\cap C^1([0,\infty); H^{s-1}(R^1))$.
\end{theorem}

Therefore  in this section we consider the case that $(u_0,v_0)\in C^{\infty}_c(R^1)$ and let $(u,v)$ be the corresponding smooth solution to (\ref{eq-dirac}) and (\ref{eq-dirac-initialv}).

Multiplying the first equation of (\ref{eq-dirac}) by $\overline{u}$   and the second equation by $\overline{v}$ gives
\begin{equation}\label{eq-dirac1}
\left\{\begin{array}{l}
(|u|^2)_t +(|u|^2)_x=mi(\overline{u}v-u\overline{v})+2\Re (i\overline{N_1}u), \\  (|v|^2)_t-(|v|^2)_x=-mi(\overline{u}v-u\overline{v})+2\Re (i\overline{N_2}v),
\end{array}\right.
\end{equation}
 which, together with the structure of nonlinear terms, leads to
\begin{equation}\label{eq-dirac-conserv} (|u|^2+|v|^2)_t +(|u|^2-|v|^2)_x=0.\end{equation}
We will use (\ref{eq-dirac-conserv}) to establish the estimates on the solution $(u,v)$ in the triangle domains as follow.

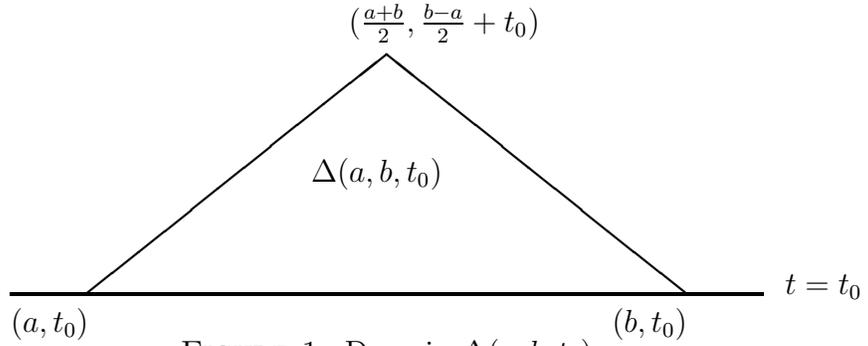
\begin{figure}[h]
\begin{center}
\unitlength=10mm
\begin{picture}(10,3.5)
\thicklines
\put(0,0){\line(1,0){10}}
\put(1,0){\line(5,4){4}}
\put(9,0){\line(-5,4){4}}
\put(0,-0.5){$(a,t_0)$}
\put(8,-0.5){$(b,t_0)$}\put(10.3,0){$t=t_0$}
\put(4.5,3.5){$(\frac{a+b}{2},\frac{b-a}{2}+t_0)$}
\put(4,1.5){$\Delta(a,b,t_0)$}
\end{picture}
\caption{Domain $\Delta(a,b,t_0)$}\label{fig-domain}
\end{center}
\end{figure}

For any $a, b\in R^1$ with $a<b$ and for any $t_0\ge 0$, we denote
\[\Delta(a,b,t_0)=\{ (x,t)\big|\, a-t_0+t<x<b+t_0-t, \, t_0<t<\frac{b-a}{2}+t_0\}, \]
 see Figure \ref{fig-domain}.
It is obvious that $\Delta(a,b,t_0)$ is bounded by two characteristic lines and $t=t_0$. The vertices of $\Delta(a,b,t_0)$ are $(a,t_0)$, $(b,t_0)$ and $(\frac{a+b}{2}, \frac{b-a}{2}+t_0)$.

\begin{lemma}\label{lemma-conv-charge-1}
For any $\tau\in [t_0, \frac{b-a}{2}+t_0]$, there holds the following,
\begin{eqnarray*}
\int\limits_{a-t_0+\tau}^{b+t_0-\tau} \big( |u(x,\tau)|^2+|v(x,\tau)|^2\big)dx +2 \int\limits_{t_0}^{\tau} |u(b+t_0-s,s)|^2 ds \\ + 2\int\limits_{t_0}^{\tau} |v(a-t_0+s,s)|^2 ds = \int\limits_a^b \big( |u(x,t_0)|^2 +|v(x,t_0)|^2\big) dx.
\end{eqnarray*}
\end{lemma}
{\it Proof.} As in  Huh \cite{huh2}, we can get the result by taking the integration  of (\ref{eq-dirac-conserv}) over the domain
\[\Omega(a,b,t_0, \tau)=\{ (x,t)\big|\, a-t_0+t<x<b+t_0-t, \, t_0<t<\tau\}. \]
The proof is complete.$\Box$

A special example of Lemma \ref{lemma-conv-charge-1} is the following estimate on the domain $\Delta(x_0-t_0, x_0+t_0,0)$ for any $x_0\in R^1$ and $t_0>0$,
\begin{eqnarray}
 2\int\limits_0^{t_0} |u(x_0+t_0-s,s)|^2ds
 +2\int\limits_0^{t_0} |v(x_0-t_0+s,s)|^2ds \nonumber\\ = \int\limits_{x_0-t_0}^{x_0+t_0} (|u_0(x)|^2+|v_0(x)|^2)dx. \label{eq-conser-charge-1}
\end{eqnarray}

\begin{lemma}\label{lemma-PointwiseEstimate}
Suppose that $\int\limits^{b}_{a}(|u_0(x)|^2+|v_0(x)|^2)dx< C_0$ for some constant $C_0>0$. Then
\begin{equation}\label{eq-pointwise1} |u(x,t)|^2 \le e^{2|\beta|C_0+mt}\big( |u_0(x-t)|^2 + mC_0\big),\end{equation}
\begin{equation}\label{eq-pointwise2} |v(x,t)|^2\le e^{2|\beta|C_0+mt}\big( |v_0(x+t)|^2 + mC_0\big),\end{equation}
for $(x,t)\in \Delta(a,b,0)$. Moreover, for any interval $[d_1,d_2]\subset [a+t,b-t]$ and $t\in [0, \frac{b-a}{2}]$, there hold the following,
\begin{equation}\label{eq-pointwise3} \int\limits_{d_1}^{d_2} |u(x,t)|^2 dx \le e^{2|\beta|C_0+mt}\big( \int\limits_{d_1-t}^{d_2-t} |u_0(x)|^2 dx +mC_0(d_2-d_1)\big),\end{equation}
\begin{equation}\label{eq-pointwise4} \int\limits_{d_1}^{d_2} |v(x,t)|^2 dx \le e^{2|\beta|C_0+mt}\big( \int\limits_{d_1+t}^{d_2+t} |v_0(x)|^2 dx +mC_0(d_2-d_1)\big).\end{equation}
\end{lemma}
 \begin{remark} We remark that the estimates (\ref{eq-pointwise1}) and (\ref{eq-pointwise2}) in Lemma \ref{lemma-PointwiseEstimate} have been proved by Huh in \cite{huh2} for $R^1\times [0,\infty)$ and we give the sketch of the proof here. Indeed, (\ref{eq-dirac1} )  gives
\begin{eqnarray*}
\frac{d}{ds} |u(x-t+s, s)|^2 &\le m\big( |u(x-t+s,s)|^2 +|v(x-t+s,s)|^2\big) \\
& +2|\beta| |u(x-t+s,s)|^2|v(x-t+s,s)|^2 ,
\end{eqnarray*}
and
\begin{eqnarray*}
\frac{d}{ds} |v(x+t-s, s)|^2 &\le m\big( |u(x+t-s,s)|^2 +|v(x+t-s,s)|^2\big) \\
& +2|\beta| |u(x+t-s,s)|^2|v(x+t-s,s)|^2.
\end{eqnarray*}
Then, taking the integration of the above from $0$ to $t$ yields that
\begin{eqnarray*}
|u(x,t)|^2 &\le& \exp( mt + 2|\beta|\int\limits_0^t|v(x-t+s,s)|^2 ds)
\\&\, &\cdot \big(|u_0(x-t)|^2 +m \int\limits_0^t e^{-ms}|v(x-t+s,s)|^2 ds \big),
\end{eqnarray*}
which implies the estimate (\ref{eq-pointwise1}) on $u$ by (\ref{eq-conser-charge-1}). The estimate (\ref{eq-pointwise2}) on $v$ could be derived in the same way. The inequalities (\ref{eq-pointwise3}) and (\ref{eq-pointwise4}) follow from (\ref{eq-pointwise1}) and (\ref{eq-pointwise2}).
\end{remark}

\subsection{Bony functional for the classical solution}

Now it follows from assumption (\ref{eq-nonlinearstruc}) on the structure of nonlinear terms that for any $(u,v)$
\begin{equation}\label{eq-nonlinearterm} |\overline{N_1}u|+ |\overline{N_2}v|\le 8|\beta||u|^2|v|^2.\end{equation}
Then, (\ref{eq-dirac1}) implies that
\begin{equation}\label{eq-dirac10}
\left\{\begin{array}{l}
(|u|^2)_t +(|u|^2)_x\le r_0(x,t), \\  (|v|^2)_t-(|v|^2)_x\le r_0(x,t).
\end{array}\right.
\end{equation}
Here
\[ r_0(x,t)= m(|u|^2+|v|^2)(x,t)+ c|u(x,t)|^2|v(x,t)|^2, \quad c=8|\beta|.\]
To deal with the  terms in the righthand sides of (\ref{eq-dirac10}),
  we define the following functionals in $\overline{\Delta (a,b,t_0)}$ as in \cite{zhang} (see also \cite{bony}, \cite{ha-tzavaras} and \cite{zhou}).
 \begin{definition}\label{def-bonyfunctional} For any $t\in [t_0, \frac{b-a}{2}+t_0]$, we define
\begin{eqnarray*}
Q_0(t,\Delta)&=&\mathop{\int\int}_{a-t_0+t<x<y<b+t_0-t} |u(x,t)|^2|v(y,t)|^2 dxdy,\\
 L_0(t,\Delta)&=&\int\limits^{b+t_0-t}_{a-t_0+t} (|u(x,t)|^2+|v(x,t)|^2)dx,\\
D_0(t,\Delta)&=&\int\limits^{b+t_0-t}_{a-t_0+t} |u(x,t)|^2|v(x,t)|^2 dx.
\end{eqnarray*}
Here $\Delta=\Delta(a,b,t_0)$.
\end{definition}
$Q_0(t,\Delta) $ is  a Bony type functional \cite{bony} (see also \cite{ha-tzavaras}, \cite{zhang} and \cite{zhou}) and is analogous to the Glimm interaction potential \cite{glimm}.
\begin{lemma}\label{lemma-bony-functional} There exists  constants $\delta_0>0$ such that for the initial data satisfying $L_0(t_0,\Delta)\le \delta_0$ there holds the following
\begin{eqnarray}\label{eq-bony-functional-ineq2}
\frac{d Q_0(t,\Delta)}{dt }+D_0(t,\Delta)\le 2m(L_0(t_0,\Delta))^2
\end{eqnarray}
for $t\in (t_0, \frac{b-a}{2}+t_0)$.
Therefore,
\begin{eqnarray}
 Q_0(t,\Delta) + \int\limits_{t_0}^t D_0(\tau,\Delta) d\tau &\le& 2m(L_0(t_0,\Delta))^2 (t-t_0) + Q_0(t_0,\Delta)  \,  \,  \,  \nonumber \\ &\le& 2m(L_0(t_0,\Delta))^2 (t-t_0)+ (L_0(t_0,\Delta))^2 \,  \,  \, \label{eq-bony-functional-ineq3}
\end{eqnarray}
for $t\in (t_0, \frac{b-a}{2}+t_0)$. Here $\Delta=\Delta(a,b,t_0)$.
\end{lemma}
{\it Proof.}  In the proof we use $L_0(t)$, $Q_0(t)$ and $D_0(t)$ for $L_0(t, \Delta)$, $Q_0(t, \Delta)$ and $D_0(t,\Delta)$ for simplification. From (\ref{eq-dirac10}) we deduce that
\begin{eqnarray*}
\frac{d Q_0(t)}{dt }  &\le& -\int\limits^{b+t_0-t}_{a-t_0+t} |u(a-t_0+t,t)|^2 |v(y,t)|^2 dy  -\int\limits^{b+t_0-t}_{a-t_0+t} |u(x,t)|^2 |v(b+t_0-t,t)|^2 dx \\
&\,&-\mathop{\int\int}_{a-t_0+t<x<y<b+t_0-t} (|u(x,t)|^2)_x |v(y,t)|^2 dxdy \\ &\,&+ \mathop{\int\int}_{a-t_0+t<x<y<b+t_0-t} |u(x,t)|^2 (|v(y,t)|^2)_y dxdy \\ & \, &+ \mathop{\int\int}_{a-t_0+t<x<y<b+t_0-t}\big(r_0(x,t)|v(y,t)|^2+ |u(x,t)|^2r_0(y,t) \big)  dxdy  \\
&\le& -2\int\limits^{b+t_0-t}_{a-t_0+t} |u(x,t)|^2|v(x,t)|^2 dx
\\
&\, & + \int\limits^{b+t_0-t}_{a-t_0+t} r_0(x,t)dx \int^{b+t_0-t}_{a-t_0+t} |v|^2 dy + \int\limits^{b+t_0-t}_{a-t_0+t} |u|^2 dx\int\limits^{b+t_0-t}_{a-t_0+t} r_0(y,t)dy \\
 &\le & \big( -2+2L_0(t)c \big)\int\limits^{b+t_0-t}_{a-t_0+t} |u(x,t)|^2|v(x,t)|^2 dx + 2m(L_0(t))^2 \\
 &= & \big( -2+2L_0(t_0)c \big)D_0(t) + 2m(L_0(t_0))^2,
\end{eqnarray*}
where  Lemma \ref{lemma-conv-charge-1}  is used in getting last equality.
Now, choose a constant $\delta_0>0$ such that
\[ -2+2\delta_0 c<-1,\]then we can get the desired result. The proof is complete. $\Box$

\subsection{Estimates on the difference between classical solutions}

 Let $(u,v)$ and $(u^{\prime}, v^{\prime})$ be the global smooth solutions to (\ref{eq-dirac}) taking $(u_0,v_0)\in C_c^{\infty}(R^1)$ and $(u^{\prime}_0, v^{\prime}_0)\in C_c^{\infty}(R^1)$ as the initial data respectively. Then the lemmas, propositions and theorems in subsection 2.1  hold for these smooth solutions $(u,v)$ and $(u^{\prime},v^{\prime})$.

Consider the difference between these two smooth solutions in the domain $\Delta(a,b,t_0)$  and let $U=u-u^{\prime}$ and $V=v-v^{\prime}$. Then,
\begin{equation}\label{eq-dirac-difference1}
\left\{ \begin{array}{l} U_t+U_x=imV -i\Big(N_1(u,v)-N_1(u^{\prime},v^{\prime})\Big), \\V_t-V_x=imU-i\Big(N_2(u,v) -N_2(u^{\prime},v^{\prime})\Big),
\end{array}
\right.
\end{equation}
which leads to
\begin{equation}\label{eq-dirac-difference2}
\left\{ \begin{array}{l} (|U|^2)_t+(|U|^2)_x=\Re 2\{imV\overline{U} -i\Big(N_1(u,v)-N_1(u^{\prime},v^{\prime})\Big)\overline{U}\}, \\(|V|^2)_t-(|V|^2)_x=\Re 2\{imU\overline{V}-i\Big(N_2(u,v) -N_2(u^{\prime},v^{\prime})\Big)\overline{V}\}
\end{array}
\right.
\end{equation}
and
\begin{equation}\label{eq-dirac-difference21}
(|U|^2+|V|^2)_t+(|U|^2-|V|^2)_x=-R
\end{equation}
with
\[ R= \Re \{2i\Big(N_1(u,v)-N_1(u^{\prime},v^{\prime})\Big)\overline{U}+ 2i\Big(N_2(u,v) -N_2(u^{\prime},v^{\prime})\Big)\overline{V}\}.\]

For the  righthand sides of (\ref{eq-dirac-difference1}-\ref{eq-dirac-difference21}), we have the following.
\begin{lemma}\label{lemma-difference}
There exists a constant $c_*>0$ such that
\[ |\Re 2\{imV\overline{U} -i\Big(N_1(u,v)-N_1(u^{\prime},v^{\prime})\Big)\overline{U}\}| \le r_1(x,t),\]
\[ |\Re 2\{imU\overline{V}-i\Big(N_2(u,v) -N_2(u^{\prime},v^{\prime})\Big)\overline{V}\}|\le  r_1(x,t),\]
and
\[ |R|\le c_*r_2(t,x,x),\]
where
\[ r_1(x,t)=m(|U(x,t)|^2+|V(x,t)|^2)+c_* r_2(t,x,x)\]
 and
\[  r_2(t,x,y)= |U(x,t)|^2\Big( |v(y,t)|^2+|v^{\prime}(y,t)|^2\Big)  + \Big(|u(x,t)|^2+|u^{\prime}(x,t)|^2\Big)|V(y,t)|^2.
\]
\end{lemma}
{\it Proof.} The assumption (\ref{eq-nonlinearstruc}) implies that there exists a constant $c_*>0$ such that
\begin{eqnarray*}
|\Big(N_1(u,v)-N_1(u^{\prime},v^{\prime})\Big)\overline{U}|
&=&|N_1(u,v)-N_1(u^{\prime},v^{\prime})||U|  \\
& \le& \frac{c_*}{8}\big(|u-u^{\prime}||v|^2+|u^{\prime}||v-v^{\prime}||v|+|u^{\prime}||v^{\prime}||v-v^{\prime}| \big)|U| \\
& \le& \frac{c_*}{4} r_2(t,x,x),
\end{eqnarray*}
and
\[ |\Big(N_2(u,v)-N_2(u^{\prime},v^{\prime})\Big)\overline{V}| \le \frac{c_*}{4} r_2(t,x,x).\]
Thus, the above estimates give the desired results. The proof is complete.$\Box$

\begin{lemma}\label{lemma-r-2}
There holds the following,
\[ |u(x,t)v(x,t)-u^{\prime}(x,t)v^{\prime}(x,t)|^2 \le 2 r_2(t,x,x)\]
for any $(x,t)\in R^1\times [0,\infty)$.
\end{lemma}
{\it Proof.} Direct computation shows that
\[ uv-u^{\prime}v^{\prime}=Uv+u^{\prime}V,\]
then the result follows. The proof is complete.$\Box$

Now to deal with the nonlinear terms $r_1$ and $r_2$, we define the following.
\begin{definition}\label{def-bony-difference}
For the domain $\Delta=\Delta(a,b,t_0)$ and for $t\in [t_0, \frac{b-a}{2}+t_0]$, define
\begin{eqnarray*} L_1(t, \Delta)&=&\int\limits^{b+t_0-t}_{a-t_0+t} (|U(x,t)|^2+|V(x,t)|^2)dx, \\
 Q_1(t,\Delta)&=&\mathop{\int\int}_{a-t_0+t<x<y<b+t_0-t} r_2(t,x,y) dxdy,\\
 D_1(t,\Delta) &=& \int\limits^{b+t_0-t}_{a-t_0+t} r_2(t,x,x) dx,
\end{eqnarray*}
where $r_2$ is given in Lemma \ref{lemma-difference}.
\end{definition}
In addition we use the notations in Definition \ref{def-bonyfunctional} for $(u,v)$, and use the following for $(u^{\prime}, v^{\prime})$,
\begin{definition}\label{def-bonyfunctional-2}
For $t\in [t_0, \frac{b-a}{2}+t_0]$, define \[ r^{\prime}_0(x,t)=m(|u^{\prime}(x,t)|^2+|v^{\prime}(x,t)|^2)+c|u^{\prime}(x,t)|^2|v^{\prime}(x,t)|^2,\]
and
\begin{eqnarray*} Q_0^{\prime}(t,\Delta)&=&\mathop{\int\int}_{a-t_0+t<x<y<b+t_0-t} |u^{\prime}(x,t)|^2|v^{\prime}(y,t)|^2 dxdy,\\
 L_0^{\prime}(t,\Delta) &=& \int\limits^{b+t_0-t}_{a-t_0+t} (|u^{\prime}(x,t)|^2+|v^{\prime}(x,t)|^2)dx,\\
 D_0^{\prime}(t,\Delta)&=&\int\limits^{b+t_0-t}_{a-t_0+t} |u^{\prime}(x,t)|^2|v^{\prime}(x,t)|^2 dx.
 \end{eqnarray*}
\end{definition}

\begin{lemma}\label{lemma-functinal-difference}
There exist constants $\delta\in (0,\delta_0)$ and $K>0$ such that if $L_0(t_0,\Delta)<\delta$ and $L_0^{\prime}(t_0,\Delta)<\delta$ then
\begin{eqnarray}
 \, & \frac{d}{dt}(L_1(t,\Delta)+KQ_1(t,\Delta))+D_1(t,\Delta) \nonumber\\ \, & \le  \Big(2mL_0(t_0,\Delta)+2mL_0^{\prime}(t_0,\Delta)+cD_0(t,\Delta)+cD^{\prime}_0(t,\Delta)\Big)L_1(t,\Delta) \label{eq-bony-functional-4}
\end{eqnarray}
for $t\in [t_0, \frac{b-a}{2}+t_0]$. Therefore
\begin{equation}\label{eq-bony-functional-4-1}
 L_1(t,\Delta)+KQ_1(t,\Delta)\le (L_1(t_0,\Delta)+KQ_1(t_0,\Delta)) \exp(h_3(t))\end{equation} and
 \begin{equation}\label{eq-bony-functional-5}
\int\limits_{t_0}^t D_1(\tau,\Delta)d\tau\le (L_1(t_0,\Delta)+KQ_1(t_0,\Delta)) \big((4m\delta+4m\delta^2c) t+2c\delta^2+1\big)\exp(h_3(t))
\end{equation}
for $t\in [t_0, \frac{b-a}{2}+t_0]$, where
\[h_3(t)=2m\big(L_0(t_0,\Delta)+L_0^{\prime}(t_0,\Delta)\big)(t-t_0) +\int\limits_{t_0}^t(cD_0(\tau,\Delta)+cD^{\prime}_0(\tau,\Delta))d\tau.\]
\end{lemma}
{\it Proof.} For simplification, in the proof we drop "$\Delta$" from the notations given  in Definitions \ref{def-bonyfunctional}, \ref{def-bony-difference} and \ref{def-bonyfunctional-2}.

 Lemma \ref{lemma-difference} and (\ref{eq-dirac-difference21}) give
\begin{equation}\label{eq-funct-linear} \frac{d}{dt}L_1(t) \le 2c_* D_1(t).\end{equation}
Next it follows from (\ref{eq-dirac-difference2}) and Lemma \ref{lemma-difference} that
\begin{equation}\label{eq-dirac-difference3}
\left\{ \begin{array}{l} (|U|^2)_t+(|U|^2)_x\le r_1(x,t), \\(|V|^2)_t-(|V|^2)_x\le r_1(x,t).
\end{array}
\right.
\end{equation}
 Then (\ref{eq-dirac-difference3}) yields that
\begin{eqnarray}
\frac{d}{dt}Q_1(t) &\le& -2D_1(t) + \mathop{\int\int}_{a-t_0+t<x<y<b+t_0-t} r_1(x,t) (|v(y,t)|^2+|v^{\prime}(y,t)|^2)dxdy  \nonumber\\
 &\,& + \mathop{\int\int}_{a-t_0+t<x<y<b+t_0-t} |U(x,t)|^2(r_0(y,t)+r_0^{\prime}(y,t))dxdy \nonumber\\
 &\,& +\mathop{\int\int}_{a-t_0+t<x<y<b+t_0-t} (|u(x,t)|^2+|u^{\prime}(x,t)|^2)r_1(y,t)dxdy  \nonumber\\
 &\,& +\mathop{\int\int}_{a-t_0+t<x<y<b+t_0-t} (r_0(x,t)+r_0^{\prime}(x,t))|V(y,t)|^2 dxdy  \nonumber\\
 &\le& -2D_1(t)+ (mL_1(t)+c_* D_1(t))(L_0(t)+L_0^{\prime}(t)) \nonumber\\
 &\,& + (mL_0(t)+mL_0^{\prime}(t)+cD_0(t)+cD^{\prime}_0(t))L_1(t) \nonumber\\
 &\le& \Big(-2+c_*(L_0(t_0)+L_0^{\prime}(t_0))\Big)D_1(t) \nonumber\\ &\,&+ \Big(2mL_0(t_0)+2mL_0^{\prime}(t_0)+cD_0(t)+cD^{\prime}_0(t)\Big)L_1(t), \label{eq-funct-quadratic}
\end{eqnarray}
where Lemma \ref{lemma-conv-charge-1} is used in getting the last inequality.

Choose small $\delta>0$ and large $K>0$ such that
\[ -2+2c_*\delta <-1\]
and \[ -K+2c_*<-1,\]
then the above estimates (\ref{eq-funct-linear}) and (\ref{eq-funct-quadratic}) give   (\ref{eq-bony-functional-4})  for $L_0(0)<\delta$ and $L_0^{\prime}(0)<\delta$. (\ref{eq-bony-functional-4-1}) is a consequence of (\ref{eq-bony-functional-4}).

To prove (\ref{eq-bony-functional-5}), we integrate (\ref{eq-bony-functional-4}) over $[t_0,t]$, then by (\ref{eq-bony-functional-4-1}) we have
\begin{eqnarray*}
\int\limits_{t_0}^t D_1(\tau)d\tau &\le& L_1(t_0)+KQ_1(t_0)+\int\limits_{t_0}^t \Big(4m\delta+cD_0(\tau)+cD^{\prime}_0(\tau)\Big)L_1(\tau)d\tau \\
 &\le& (L_1(t_0)+KQ_1(t_0)) \\ &+ &\int\limits_{t_0}^t\Big(4m\delta+cD_0(\tau)+cD^{\prime}_0(\tau)\Big) d\tau (L_1(t_0)+KQ_1(t_0))\exp (h_3(t))\\
 &\le& (L_1(t_0)+KQ_1(t_0)) \big((4m\delta+4m\delta^2c) t+2c\delta^2+1\big)\exp(h_3(t)),
\end{eqnarray*}
where the estimates (\ref{eq-bony-functional-ineq3}) for $D_0$ and $D_0^{\prime}$ are also used for the last inequality.
 This completes the proof. $\Box$

\begin{lemma}\label{lemma-h3(t)}
Let $h_3(t)$ be given by Lemma \ref{lemma-functinal-difference}. Then for $t\ge t_0$,
\begin{eqnarray*}
 h_3(t)\le 4m(\delta+\delta^2)(t-t_0) +2\delta^2.
\end{eqnarray*}
\end{lemma}
{\it Proof.} We use the notations in the proof of Lemma \ref{lemma-functinal-difference} for simplification.

 By Lemma \ref{lemma-bony-functional}, we have
\begin{eqnarray*}
h_3(t)& \le & 2m\Big(L_0(t_0)+L_0^{\prime}(t_0)+(L_0(t_0))^2+ (L_0^{\prime}(t_0))^2 \Big)(t-t_0)
\\ &+& (L_0(t_0))^2+(L_0^{\prime}(t_0))^2,
\end{eqnarray*}
which yields the desired estimate for $h_3$. The proof is complete. $\Box$

  \section{Proof of Theorem \ref{thm-existence}}\label{section-proof}

  \subsection{Existence of solution}

   Choose a sequence of smooth functions \[(u_{k0},v_{k0})\in C_c^{\infty}(R^1), \quad k=1,2\cdots,\] such that
\begin{equation}\label{eq-initial-approx} \lim\limits_{k\to\infty}(||u_{k0}-u_{0}||_{L^2(R^1)}^2 +||v_{k0}-v_{0}||_{L^2(R^1)}^2)=0.\end{equation}
Theorem \ref{thm-huh} implies that there is a sequence of smooth solutions, $(u_k,v_k)$, $k=1,2\cdots$, which take $(u_{k0},v_{k0})\in C_c^{\infty}(R^1),$ $k=1,2\cdots,$ as their initial data respectively. Moreover $(u_k(\cdot, t), v_k(\cdot,t))$ as function of $x$ has compact support for each $t\ge 0$ and for $k\ge 1$.

We will show that the sequence of solutions $(u_k,v_k)$ converges to a strong solution in the domain $R^1\times [0,T]$ for any $T>0$. To this end, we first choose a $B>0$ such that the following hold,
\[ \int\limits_{|x|\ge B}( |u_0(x)|^2+|v_0(x)|^2) dx <\frac{\delta}{4}.\]
Moreover, due to the convergence of the sequence $\{(u_{k0},v_{k0}) \}_{k=1}^{\infty}$,  $B$ can be chosen so that
\[\sup_{k\ge 1} \int\limits_{|x|\ge B}( |u_{k0}(x)|^2+|v_{k0}(x)|^2) dx <\frac{\delta}{3}.\]
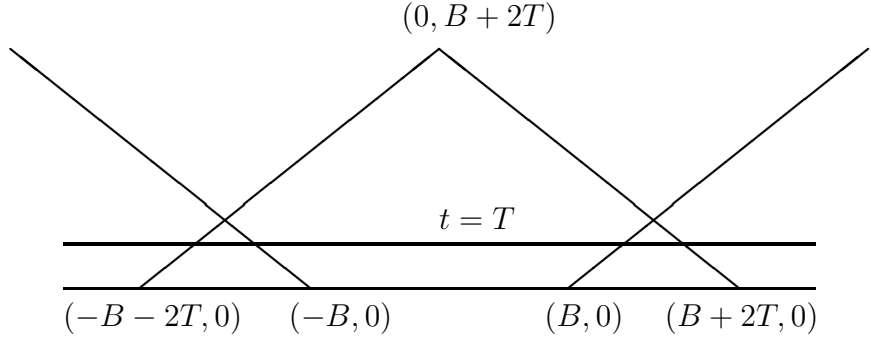
\begin{figure}[h]
\begin{center}
\unitlength=10mm
\begin{picture}(12,6.5)
\thicklines
\put(0,2){\line(1,0){10}}
\put(0,2.6){\line(1,0){10}}
\put(1,2){\line(5,4){4}}
\put(3.3,2){\line(-5,4){4}}\put(3,1.5){$(-B,0)$}
\put(9,2){\line(-5,4){4}}
\put(6.7,2){\line(5,4){4}}\put(6.4,1.5){$(B,0)$}
\put(0,1.5){$(-B-2T,0)$}
\put(8,1.5){$(B+2T,0)$}
\put(5,2.8){$t=T$}
\put(4.5,5.5){$(0,B+2T)$}
\end{picture}
\caption{Division of domain}\label{fig-Division-of-domain}
\end{center}
\end{figure}

\begin{lemma}\label{lemma-convergence1}
For any $T>0$, there hold the following,
\[ \lim_{k,l\to\infty} \max_{t\in [0,T]} \int\limits_{|x|\ge B+t} \big(|u_k(x,t)-u_l(x,t)|^2+ |v_k(x,t)-v_l(x,t)|^2\big)dx = 0\]
and
\[ \lim_{k,l\to\infty} \mathop{\int\int}_{t\in [0,T], \, \, |x|\ge B+t} |u_kv_k(x,t)-u_lv_l(x,t)|^2dxdt = 0.\]
\end{lemma}
{\it Proof.} For any $k,l$,  $u_k$ and $u_l$ have compact supports in $R^1\times [0,T]$. Then
\[ supp_T(u_k)\cup supp_T(u_l)\subset \Big(\Delta(-B_{k,l}, B_{k,l},0)\cap (R^1\times [0,T])\Big)\]
for some $B_{k,l}>0$, here
\[supp_T (u_k)= \overline{\{ (x,t)|\, x\in R^1, \, t\in [0,T],\, u_k(x,t)>0\}} \] and \[ supp_T (u_l)= \overline{\{ (x,t)|\, x\in R^1, \, t\in [0,T],\, u_l(x,t)>0\}}.\]
Due to the choice of $B$, we have
\[ \int\limits_{-B_{k,l}}^{-B} \big(|u_{k0}(x)|^2+|u_{l0}(x)|^2+|v_{k0}(x)|^2+|v_{l0}(x)|^2 \big) dx<\delta,\]
and
\[ \int\limits^{B_{k,l}}_{B} \big(|u_{k0}(x)|^2+|u_{l0}(x)|^2+|v_{k0}(x)|^2+|v_{l0}(x)|^2 \big) dx<\delta,\]
Now, applying Lemma \ref{lemma-r-2} and Lemma \ref{lemma-functinal-difference} to the functions $(u_k,v_k)$ and $(u_l,v_l)$ on the domains $\Delta(-B_{k,l}, -B,0)$ and $\Delta(B,B_{k,l},0)$ yields that
\begin{eqnarray*}
 \max_{t\in [0,T]} \int\limits_{|x|\ge B+t} \big(|u_k(x,t)-u_l(x,t)|^2+ |v_k(x,t)-v_l(x,t)|^2\big)dx \\
 +
 \mathop{\int\int}_{t\in [0,T], \, \, |x|\ge B+t} |u_kv_k(x,t)-u_lv_l(x,t)|^2dxdt \\ \le C(T) \int\limits_{|x|\ge B} \big(|u_{k0}(x)-u_{l0}(x)|^2+|v_{k0}(x)-v_{l0}(x)|^2 \big) dx,
 \end{eqnarray*}
 where the constant $C(T)$ can be determined due to Lemmas \ref{lemma-bony-functional}, \ref{lemma-functinal-difference} and Lemma \ref{lemma-h3(t)} and is independent of $k,l$. This gives the desired result and completes the proof. $\Box$

Now we  want to prove the following.
\begin{proposition}\label{prop-convergence2}
For any $T>0$, there hold the following,
\[ \lim_{k,l\to \infty}\max_{t\in[0,T]} ||u_k(\cdot, t)-u_l(\cdot,t)||_{L^2(R^1)}=0,\]
and
\[ \lim_{k,l\to \infty}||u_kv_k-u_lv_l||_{L^2(R^1\times [0,T])} =0.\]
\end{proposition}
{\it Proof.} By Lemma \ref{lemma-convergence1}, it suffices to prove the convergence of the sequences $\{u_k\}$ and $\{u_kv_k\}$ in the domain $\Delta(-B-2T, B+2T, 0)$. We first choose the constant $C_0$ in Lemma \ref{lemma-PointwiseEstimate} as
\[ C_0=1+\sup_{k\ge 1} \int\limits^{\infty}_{-\infty} (|u_{k0}|^2+|v_{k0}|^2)dx +\int\limits_{-\infty}^{\infty}(|u_0|^2+|v_0|^2)dx.\]
Moreover, there exists a $r>0$ such that for any interval $[a_1,a_2]\subset [-B-4T, B+4T]$ with the length $|a_2-a_1|\le 4r$ there holds that
\begin{equation}\label{eq-small-interval1} e^{2|\beta|C_0+mT}\big(  \int\limits_{a_1}^{a_2}(|u_0(x)|^2+|v_0(x)|^2)dx+mC_0(a_2-a_1)\big)<\frac{\delta}{8}.\end{equation}
Due to the convergence of the sequence $\{ (u_{k0},v_{k0})\}_{k=1}^{\infty}$, $r$ could be chosen so that
\begin{equation}\label{eq-small-interval2} e^{2|\beta|C_0+mT}\big(  \int\limits_{a_1}^{a_2}(|u_{k0}(x)|^2+|v_{k0}(x)|^2)dx+mC_0(a_2-a_1)\big)<\frac{\delta}{4},\end{equation}
and
 $B+2T=2Nr$ for some positive integer $N$.

\begin{figure}[h]
\begin{center}
\unitlength=10mm
\begin{picture}(12,6.5)
\thicklines
\put(0,2){\line(1,0){10}}
\put(0,3.7){\line(1,0){10}}
\put(1,2){\line(5,4){4}}
\put(9,2){\line(-5,4){4}}
\put(0,1.5){$(mr,jr)$}
\put(8,1.5){$((m+4)r,jr)$}
\put(0,4){$t=(j+1)r$}
\put(3.1,0.5){\line(0,1){4.5}}\put(2,0){$x=(m+1)r$}
\put(6.9,0.5){\line(0,1){4.5}}\put(5.8,0){$x=(m+3)r$}
\put(4.5,5.5){$((m+2)r,(j+2)r)$}
\end{picture}
\caption{Domain $\Delta(mr,(m+4)r, jr)$}\label{fig-induction}
\end{center}
\end{figure}
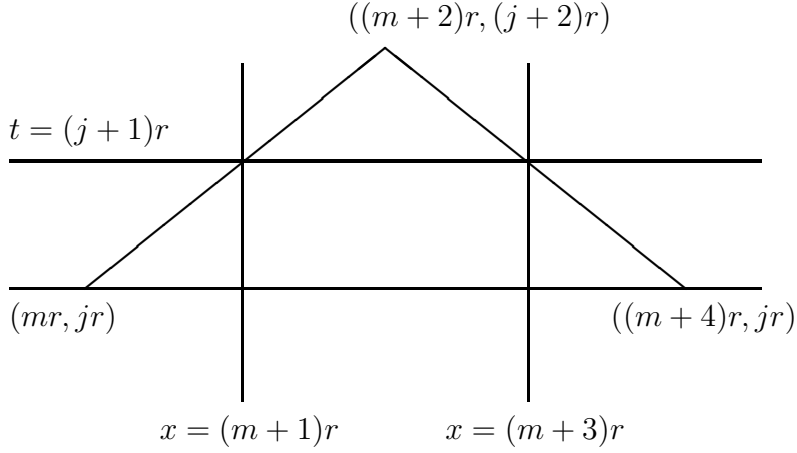

Choose grid points $P_{m,j}=(mr, jr)$, $-2N\le m\le 2N$, $0\le j\le 2N$, in the domain $\Delta(-B-2T,B+2T,0)$, see Figure \ref{fig-induction}.

Then, with (\ref{eq-initial-approx}), (\ref{eq-small-interval1}) and (\ref{eq-small-interval2}),  applying  Lemmas \ref{lemma-functinal-difference} and \ref{lemma-h3(t)} to $(u_k,v_k)$ and $(u_l,v_l)$ in the domain $\Delta(mr,(m+4)r, 0)$ for $k\ge 1$ and $-2N\le m\le 2N-4$ yields that
\begin{eqnarray}
\max_{0\le t\le 2r}\int\limits_{mr+t}^{(m+4)r-t}\big(|u_k(x, t)-u_l(x,t)|^2+|v_k(x, t)-v_l(x,t)|^2\big) dx \nonumber \\
+\mathop{\int\int}_{0\le t\le 2r, \,  \, \,  mr+t\le x \le (m+4)r-t} |u_kv_k(x,t)-u_lv_l(x,t)|^2 dxdt  \nonumber \\
\le C(T)\int\limits_{mr}^{(m+4)r} \big(|u_k(x,0)-u_l(x,0)|^2+|v_k(x, 0)-v_l(x,0)|^2\big) dx \to 0 \label{eq-convergence-0m}
\end{eqnarray}
as $k,l\to \infty$. Here the constant  $C(T)$ depends only on $T$ and is determined by Lemmas \ref{lemma-bony-functional},  \ref{lemma-functinal-difference} and \ref{lemma-h3(t)}.

Taking the summation of (\ref{eq-convergence-0m}) over $m$ yields that
\[\displaystyle \max_{0\le t\le r} \int\limits_{-B-2T+t}^{B+2T-t}\big(|u_k(x, t)-u_l(x,t)|^2+|v_k(x, t)-v_l(x,t)|^2\big) dx \to 0\]
and
\[\mathop{\int\int}_{0\le t\le r, \, \, \,  -B-2T+t\le x \le B+2T-t} |u_kv_k(x,t)-u_lv_l(x,t)|^2 dxdt \to 0
\] as $k,l\to \infty$.

Now the proof can be carried out by induction on $jr$ by using the set of domains $\{\Delta(mr,(m+4)r, jr);-2N+j\le m\le 2N-4-j\}$. Indeed suppose that
\[ \lim_{k,l\to \infty}\max_{0\le t \le jr}\int\limits_{-B-2T+t}^{B+2T-t} \big(|u_k(x, t)-u_l(x,t)|^2+|v_k(x, t)-v_l(x,t)|^2\big) dx=0\]
and
\[ \lim_{k,l\to\infty} \mathop{\int\int}_{0\le t\le jr, \, \, \, -B-2T+t<x<B+2T-t} |u_kv_k(x,t)-u_lv_l(x,t)|^2 dxdt=0.\]
By Lemma \ref{lemma-PointwiseEstimate} and by (\ref{eq-small-interval2}) we have
\[ \int\limits_{mr}^{(m+4)r} |u_k(x,jr)|^2 dx <\frac{\delta}{2}, \quad \int\limits_{mr}^{(m+4)r} |v_k(x,jr)|^2 dx <\frac{\delta}{2} \]for $k\ge 1, \, \, -2N+j\le m\le 2N-4-j.$

Then, applying Lemma \ref{lemma-functinal-difference} to $(u_k,v_k)$ and $(u_l,v_l)$ in the domain $\Delta(mr,(m+4)r, jr)$ for $k\ge 1$ and $-2N+j\le m\le 2N-4-j$ yields that
\begin{eqnarray}
\max_{jr\le t\le (j+2)r}\int\limits_{mr-jr+t}^{(m+4)r+jr-t}\big(|u_k(x, t)-u_l(x,t)|^2+|v_k(x, t)-v_l(x,t)|^2\big) dx \nonumber \\
+\mathop{\int\int}_{jr\le t\le (j+2)r, \,  \, \,  mr-jr+t\le x \le (m+4)r+jr-t} |u_kv_k(x,t)-u_lv_l(x,t)|^2 dxdt  \nonumber \\
\le C(T)\int\limits_{mr}^{(m+4)r} \big(|u_k(x, jr)-u_l(x,jr)|^2+|v_k(x, jr)-v_l(x,jr)|^2\big) dx \to 0 \label{eq-convergence-jm}
\end{eqnarray}
as $k,l\to \infty$. Here the constant  $C(T)$ depends only on $T$ and is determined by Lemmas \ref{lemma-bony-functional},  \ref{lemma-functinal-difference} and \ref{lemma-h3(t)}.

Then taking the summation of (\ref{eq-convergence-jm}) over $m$ yields that
\[\displaystyle \max_{jr\le t\le (j+1)r} \int\limits_{-B-2T+t}^{B+2T-t}\big(|u_k(x, t)-u_l(x,t)|^2+|v_k(x, t)-v_l(x,t)|^2\big) dx \to 0\]
and
\[\mathop{\int\int}_{jr\le t\le (j+1)r, \, \, \,  -B-2T+t\le x \le B+2T-t} |u_kv_k(x,t)-u_lv_l(x,t)|^2 dxdt \to 0
\] as $k,l\to \infty$. Therefore, by induction hypothesis, it leads to the convergence in $L^2(\Delta(-B-2T,B+2T,0))$ of the sequences $\{u_k\}$ and $\{u_kv_k\}$. The proof is complete. $\Box$

\subsection{Uniqueness of solution and Proof of main result}

\begin{proposition}\label{prop-unique}
Suppose that $\{(u^{\prime}_k,v^{\prime}_k)\}$ and $\{(u^{\prime\prime}_k,v^{\prime\prime}_k)\}$ are two sequences of smooth solutions to (\ref{eq-dirac}) with
\[ \lim_{k\to \infty}\int\limits^B_{-B} \big(|u_{k0}^{\prime}(x)-u^{\prime\prime}_{k0}(x)|^2+|v_{k0}^{\prime}(x)-v^{\prime\prime}_{k0}(x)|^2 \big)dx =0 \] for some $B>0$.  Here $(u_{k0}^{\prime},v_{k0}^{\prime})$ and $(u^{\prime\prime}_{k0},v^{\prime\prime}_{k0})$ are initial data of $(u^{\prime}_k,v^{\prime}_k)$ and $(u^{\prime\prime}_k,v^{\prime\prime}_k)$ respectively. Then
\[ \lim_{k\to\infty} \max_{0\le t\le B} \int\limits_{-B+t}^{B-t} \big(|u_{k}^{\prime}(x,t)-u^{\prime\prime}_{k}(x,t)|^2+|v_{k}^{\prime}(x,t)-v^{\prime\prime}_{k}(x,t)|^2 \big)dx =0. \]
\end{proposition}
{\it Proof.} By the characteristic method, we can show that the smooth solutions $\{(u^{\prime}_k,v^{\prime}_k)\}$ and $\{(u^{\prime\prime}_k,v^{\prime\prime}_k)\}$ in $\Delta(-B,B,0)$ are uniquely determined by the data $(u_{k0}^{\prime}, v)_{k0}^{\prime})\big|_{[-B,B]}$ and $(u_{k0}^{\prime\prime}, v)_{k0}^{\prime\prime})\big|_{[-B,B]}$. Therefore, we may assume in the proof that $(u_{k0}^{\prime}, v)_{k0}^{\prime})\big|_{|x|\ge 2B}=(u_{k0}^{\prime\prime}, v)_{k0}^{\prime\prime})\big|_{|x|\ge 2B}=0$, and the lemmas in Section \ref{section-estimates-solu} could be applied to these solutions.

We choose the constant $C_0$ in Lemma \ref{lemma-PointwiseEstimate}  as
\[ C_0=1+\sup_{k\ge 1}\int\limits_{-B}^{B}\big( |u_{0k}^{\prime}|^2+|v_{0k}^{\prime}|^2+|u_{0k}^{\prime\prime}|^2+|v_{0k}^{\prime\prime}|^2\big)dx.\]
Moreover by the convergence of the sequence $\{ (u_{k0}^{\prime},v_{k0}^{\prime})\}$ and $\{ (u_{k0}^{\prime\prime},v_{k0}^{\prime\prime})\}$, there exists a $r>0$ such that for any interval $[a_1,a_2]\subset [-B, B]$ with the length $|a_2-a_1|\le 4r$ there holds that
\begin{equation}\label{eq-small-interval4} e^{2|\beta|C_0+mT}\big(  \int\limits_{a_1}^{a_2}(|u_{k0}^{\prime}(x)|^2+|v_{k0}^{\prime}(x)|^2)dx+mC_0(a_2-a_1)\big)<\frac{\delta}{8},\end{equation}
\begin{equation}\label{eq-small-interval5} e^{2|\beta|C_0+mT}\big(  \int\limits_{a_1}^{a_2}(|u_{k0}^{\prime\prime}(x)|^2+|v_{k0}^{\prime\prime}(x)|^2)dx+mC_0(a_2-a_1)\big)<\frac{\delta}{8}.\end{equation}
 In addition the length $r$ could be chosen so that
 $B=2Nr$ for some positive integer $N$.

Then, by (\ref{eq-small-interval4}) and (\ref{eq-small-interval5}), we apply Lemmas \ref{lemma-functinal-difference} and \ref{lemma-h3(t)} to $(u_k^{\prime},v_k^{\prime})$ and $(u_k^{\prime\prime},v_k^{\prime\prime})$ in the domain $\Delta(mr,(m+r)r,0)$ $(-2N\le m\le 2N-4)$ to get the following,
\begin{eqnarray*}
\max_{0\le t\le 2r}\int\limits_{mr+t}^{(m+4)r-t}\big(|u_k^{\prime}(x, t)-u_k^{\prime\prime}(x,t)|^2+|v_k^{\prime}(x, t)-v_k^{\prime\prime}(x,t)|^2\big) dx \\
\le C(T)\int\limits_{mr}^{(m+4)r} \big(|u_k^{\prime}(x, 0)-u_k^{\prime\prime}(x,0)|^2+|v_k^{\prime}(x, 0)-v_k^{\prime\prime}(x,0)|^2\big) dx \to 0
\end{eqnarray*}
as $k\to \infty$, which yields that
\[ \lim_{k\to\infty} \max_{0\le t\le r} \int\limits_{-B+t}^{B-t} \big(|u_{k}^{\prime}(x,t)-u^{\prime\prime}_{k}(x,t)|^2+|v_{k}^{\prime}(x,t)-v^{\prime\prime}_{k}(x,t)|^2 \big)dx =0. \]

Now the proof can be carried out by induction on $jr$ by using the set of domains $\{\Delta(mr,(m+4)r, jr); -2N+j\le m \le 2N-4-j\}$. Indeed suppose that
\[ \lim_{k\to \infty}\max_{0\le t\le jr}\int\limits_{-B+t}^{B-t} \big(|u_k^{\prime}(x, t)-u_k^{\prime\prime}(x,t)|^2+|v_k^{\prime}(x, t)-v_k^{\prime\prime}(x,t)|^2\big) dx=0.\]
By Lemma \ref{lemma-PointwiseEstimate} and  by (\ref{eq-small-interval4}) and (\ref{eq-small-interval5}), we have
\[ \int\limits_{mr}^{(m+4)r} |u_k^{\prime}(x,jr)|^2 dx <\frac{\delta}{2}, \quad \int\limits_{mr}^{(m+4)r} |v_k^{\prime}(x,jr)|^2 dx <\frac{\delta}{2} \]
and
\[ \int\limits_{mr}^{(m+4)r} |u_k^{\prime\prime}(x,jr)|^2 dx <\frac{\delta}{2}, \quad \int\limits_{mr}^{(m+4)r} |v_k^{\prime\prime}(x,jr)|^2 dx <\frac{\delta}{2} \]
for $k\ge 1, \, \, -2N+j\le m\le 2N-4-j.$

Then, applying Lemma \ref{lemma-functinal-difference} to $(u^{\prime}_k,v^{\prime}_k)$ and $(u^{\prime\prime}_k,v^{\prime\prime}_k)$ in the domain $\{\Delta(mr,(m+4)r, jr)\}$ for $k\ge 1$ and $-2N+j\le m\le 2N-4-j$ yields that
\begin{eqnarray*}
\max_{jr\le t\le (j+2)r}\int\limits_{mr-jr+t}^{(m+4)r+jr-t}\big(|u_k^{\prime}(x, t)-u_k^{\prime\prime}(x,t)|^2+|v_k^{\prime}(x, t)-v_k^{\prime\prime}(x,t)|^2\big) dx \\
\le C(T)\int\limits_{mr}^{(m+4)r} \big(|u_k^{\prime}(x, jr)-u_k^{\prime\prime}(x,jr)|^2+|v_k^{\prime}(x, jr)-v_k^{\prime\prime}(x,jr)|^2\big) dx \to 0
\end{eqnarray*}
as $k\to \infty$. Here the constant  $C(T)$ is  determined by Lemmas \ref{lemma-bony-functional},  \ref{lemma-functinal-difference} and \ref{lemma-h3(t)}.

By the induction hypothesis it follows that
\[ \lim_{k\to\infty} \max_{0\le t\le (j+1)r} \int\limits_{-B+t}^{B-t} \big(|u_{k}^{\prime}(x,t)-u^{\prime\prime}_{k}(x,t)|^2+|v_{k}^{\prime}(x,t)-v^{\prime\prime}_{k}(x,t)|^2 \big)dx =0. \]
 Therefore we prove the result.  The proof is complete. $\Box$

{\bf Proof of Theorem \ref{thm-existence}:} The existence of strong solution $(u,v)$ follows from Propositions \ref{prop-convergence2}.

 For the uniqueness, let $(u^{\prime},v^{\prime})$ and $(u^{\prime\prime},v^{\prime\prime})$ be two strong solutions to (\ref{eq-dirac}) and (\ref{eq-dirac-initialv}), our aim is to prove $(u^{\prime},v^{\prime})=(u^{\prime\prime},v^{\prime\prime})$. To this end, let  $\{(u^{\prime}_k,v^{\prime}_k)\}$ and $\{(u^{\prime\prime}_k,v^{\prime\prime}_k)\}$ be two sequences of smooth solutions to (\ref{eq-dirac}), which are  convergent in $L^2(R^1\times [0,T])$ to $(u^{\prime},v^{\prime})$ and $(u^{\prime\prime},v^{\prime\prime})$ respectively for any $T>0$. Moreover their initial data  $(u_{k0}^{\prime},v_{k0}^{\prime})$ and $(u^{\prime\prime}_{k0},v^{\prime\prime}_{k0})$ are also assumed to be convergent to same data $(u_0,v_0)$.

Then by Proposition \ref{prop-unique},
\[ \lim_{k\to\infty} \mathop{\int\int}_{\Delta(-B,B,0)}\big( |u_k^{\prime}(x,t)-u_k^{\prime\prime}(x,t)|^2 +|v_k^{\prime}(x,t)-v_k^{\prime\prime}(x,t)|^2\big) dxdt =0\] for any $B>0$,
which yields that
\[ (u^{\prime},v^{\prime})(x,t)=(u^{\prime\prime},v^{\prime\prime})(x,t), \quad \mbox{for $(x,t)\in\Delta(-B,B,0)$.}\]
This leads to the uniqueness.
The proof is complete.$\Box$

\section*{ Acknowledgement}

This work was partially  supported by NSFC Project 11031001 and 11121101, by the 111 Project
B08018 and by Ministry of Education of China.

\end{document}